\documentclass[12pt,reqno]{amsart}
\usepackage{amssymb,amsmath,amsthm,longtable,multirow,array}

\usepackage{color}
\usepackage{mathrsfs}
\numberwithin{equation}{section}

\newtheorem{Theorem}{Theorem}[section]
\newtheorem{Proposition}[Theorem]{Proposition}

\newtheorem{Lemma}[Theorem]{Lemma}
\newtheorem{Corollary}[Theorem]{Corollary}
\theoremstyle{definition}
\newtheorem{Definition}[Theorem]{Definition}
\theoremstyle{remark}
\newtheorem{Example}[Theorem]{Example}

\newtheorem*{remark*}{Remark}

 \providecommand\CC{{\mathbb C}}

 \providecommand\QQ{{\mathbb Q}}
\providecommand\RR{{\mathbb R}}
\providecommand\TT{{\mathbb T}}
 \providecommand\ZZ{{\mathbb Z}}

\providecommand\brs{\begin{remark*}}
\providecommand\ers{\end{remark*}}

\providecommand\om{\omega}

\providecommand\lm{\lambda}

\providecommand\be{\begin{enumerate}}
\providecommand\ee{\end{enumerate}}
\providecommand\bT{\begin{Theorem}}
\providecommand\eT{\end{Theorem}}
\providecommand\bP{\begin{Proposition}}
\providecommand\eP{\end{Proposition}}
\providecommand\bD{\begin{Definition}}
\providecommand\eD{\end{Definition}}
\providecommand\bE{\begin{Example}}
\providecommand\eE{\end{Example}}
 \providecommand\bL{\begin{Lemma}}
\providecommand\eL{\end{Lemma}}
\providecommand\bC{\begin{Corollary}}
\providecommand\eC{\end{Corollary}}

\providecommand\bpp{\begin{proof}} \providecommand\epp{\end{proof}}
\providecommand\bee{\begin{equation}}
\providecommand\eee{\end{equation}}

  \providecommand\beqq{\begin{eqnarray*}}
\providecommand\eeqq{\end{eqnarray*}}
\providecommand\a{\alpha}

\providecommand\th{\theta}
 
\providecommand\ol{\overline}
 \providecommand\rt{\rightarrow}
\providecommand\?{\infty}

\providecommand\besp{\begin{split}}
\providecommand\eesp{\end{split}}
\providecommand\bay{\begin{array}}
\providecommand\eay{\end{array}}

\providecommand\bes{\begin{equation}\begin{split}}
\def\a{\alpha}
\def\b{\beta}

\def\f#1#2{\frac{#1}{#2}}
\def\t{\tau}
\def\th{\theta}
\def\s{\sigma}

\def\<{\leq}
\begin{document}

\title[Hilbert Transform]{Hilbert Transformation and Representation of  $ax+b$ Group}
\author{Pei Dang}
\author{Hua Liu}
\author{Tao Qian}

\begin{abstract}
In this paper we study the Hilbert transformations over $L^2(\RR)$ and $L^2(\TT)$  from
the viewpoint of symmetry. For a linear operator over $L^2(\RR)$  commutative with the ax+b group we show that the operator is of the form
$
\lambda I+\eta H,
$
where $I$ and $H$ are the identity operator and Hilbert transformation respectively, and $\lambda,\eta$ are complex numbers. In the related literature this result was proved through first invoking the boundedness result of the operator, proved though a big machinery. In our setting the boundedness is a consequence of the boundedness of the Hilbert transformation.  The methodology that we use is Gelfand-Naimark's representation of the ax+b group. Furthermore we prove a similar result on the unit circle. Although there does not exist a group like ax+b on the unit circle,  we construct a semigroup to play the same symmetry role for the Hilbert transformations over the circle $L^2(\TT).$
\end{abstract}

\address{pdang@must.edu.mo,  Faculty of Information Technology, Macau University of Science and Technology,
Macau }
\address{daliuhua@163.com, Department of Mathematics, Tianjin
University of Technology and Education, Tianjin 300222, China}
\address{fsttq@umac.mo, Department of Mathematics, Macau University}
\address{This work was supported in part by Macao Science and Technology Development Fund, MSAR. Ref. 045/2015/A2; NSFC grant 11471250; Macao Government FDCT 098/2012/A3;
University of Macau Multi-Year Research Grant (MYRG) MYRG116(Y1-L3)-FST13-QT.}

 \maketitle

\bigskip
\section{ Introduction}
\noindent
 The Hilbert transformation given by the formula
\begin{equation}\label{def-hilbert}
  \mathbf{H}f(x)=\f1{\pi}\text{p.v.}\int_\RR\f{f(y)}{x-y}dy,\quad x\in \RR,
\end{equation}
can be first defined for functions $f$ of finite energy and locally of the H\"older type continuity. It then can be extended to become a $L^2$-bounded linear operator over the whole $L^2$ space.
In the rest of the article we will omit the prix p.v. in the case of no confusion.

The Hilbert transformation has many applications, including solving problems in aerodynamics, condensed matter physics, optics, fluids, and engineering (see, for instance, \cite{king}). It is, especially, an indispensable tool in harmonic and signal analysis.

 Hilbert transformations play a role to connect harmonic with complex analysis. In general terms, a Hilbert transformation on a manifold, can be defined as the mapping from the scalar part (real part) to the non-scalar part (imaginary part) of the boundary limits of complex analytic functions on one of the two regions divided by the manifold (\cite{be}). For example, $\forall u\in L^2(\RR), u+i\mathbf{H}u$ belongs to the closed subspace of $L^2(\RR)$ constituted by the non-tangential boundary limits of the functions in the complex Hardy space $H^2(\CC^+).$  The concerned closed subspace of $L^2(\RR)$ is denoted by $H^2_+(\RR ).$ On the contrary, if $ f\in H^2(\CC^+),$ then there exists a function $u\in L^2(\RR)$ such that the non-tangential boundary limit of $f,$ still denoted as $f,$ possesses the form $u+i\mathbf{H}u,$ where $u$ can be chosen as real-valued or complex valued. In particular,  the non-tangential  boundary limit of a Hardy $H^2(\CC^+)$ function $f$ can have the expression $(1/2)f+i(1/2)\mathbf{H}f$ on the boundary, phrased as the Plemelj formula \cite{miklin}, that further implies $\mathbf{H}f=-if.$ The latter turns to be a characterization of a function $f\in L^2(\RR )$ to be in $H^2_+(\RR )$ (\cite{Ga}). The above relations can all be extended to Hardy spaces $H^p(\CC^+)$ with $1\leq p<\infty$ (\cite{qian}). For $p=2$ the trace operator taking $H^2(\CC^+)$ to $H^2_+(\RR )$ is, in fact, an isometry between the Hilbert spaces. In harmonic and signal analysis there is correspondence between the real signal $u$ and analytic one $u+i\mathbf{H}u$ (\cite{qian0, hwy}). In the current study we restrict ourselves to $p=2.$ For the lower half complex plane we have an analogous theory and, correspondingly, have the spaces $H^2(\CC^-)$ and $H^2_-(\RR).$

 Denote by $H^\pm=H^2_\pm (\RR )$ the spaces consisting of, respectively, the non-tangential boundary limits of the upper and lower half Hardy spaces, the Hilbert transformation
$\mathbf{H}$ can be decomposed into the sum of the two projection operators over $H^\pm,$ respectively. In fact, by the Plemelj formula, $$L^2(\RR)=H^+\bigoplus H^-,  \quad \mathbf{H}=(-i)(\mathbf{P}^+-\mathbf{P}^-),$$
 where $P^\pm$ denote the projection operators over $H^\pm$ respectively. By this decomposition it is easy to check that the Hilbert transformation is a power  self-inverse, and precisely,  $$(i\mathbf{H})^2=\mathbf{I}, \quad {\rm or}\quad  \mathbf{H}^2=-\mathbf{I},$$ where $\mathbf{I}$ is the identical operator \cite{miklin}.
We further note that
\[ \mathbf{H}f(x)=\f1{2\pi}\int_\RR e^{i\xi x}(-i{\rm sgn} (\xi))f^\wedge(\xi)d\xi,\]
where ${f}^\wedge$ is the Fourier transform of $f,$  and, for almost all $x,$
\begin{eqnarray*} \mathbf{P}^\pm f(x)&=&\pm \frac{1}{2\pi} \int_0^{\pm\infty} e^{i\xi x}{f}^\wedge(\xi)d\xi\\
&=&\frac{1}{2}f(x)\pm i\frac{1}{2}\mathbf{H}f(x).\end{eqnarray*}

In \cite{va} a set of characterization conditions for an operator to be the Hilbert transformation is given, while the characterization conditions are in terms of properties of the images of the operator restricted to the exponential functions. In \cite{hwy} the authors further study the aspect and give mathematical proofs of the results in \cite{va}. In the present study we give characterizations of the Hilbert transformations on the real line and on the unit circle in terms of group symmetry in the respective contexts. What is interesting is that the Hilbert transformation operator originally defined as an analysis object can be fully characterized through algebraic operations.
\vspace{2mm}

On the real line we consider the following two group operations in relation to symmetry properties of operators.
\vspace{2mm}

Denote by $T_a$ the dilation operator
\begin{equation}\label{dilation}
  T_af=a^{-\f12}f(\f x a),\ a>0,\quad \forall f\in L_2(\RR),
\end{equation}
and $\tau_b$ the translation operator
\begin{equation}\label{translation}
  \tau_bf=f(x-b),\ b\in \RR,\quad \forall f\in L_2(\RR).
\end{equation}
It is evident that
both $T_a$ and $\tau_b$ are isometric mappings from $L^2(R)$ to itself.
 Then we have the following lemma, \cite{king}.
 \bL\label{commute} For  $a\in \RR^+, b\in \RR$, both operators $T_a$ and $\tau_b$ commutes with the Hilbert
transformation $\mathbf{H},$
\begin{equation}\label{dilation-hil}
  T_a\mathbf{H}=\mathbf{H}T_a
\end{equation}
and
\begin{equation}\label{tran-hil}
  \tau_b\mathbf{H}=\mathbf{H}\tau_b.
\end{equation}
\eL

 Lemma \ref{commute}  reveals  the physical  significance of the Hilbert transformation. (\ref{dilation-hil}) means that $\mathbf{H}$ is independent of the scales, while (\ref{tran-hil}) says  that  $\mathbf{H}$  is independent of the location of the original point. We know that translation and dilation generate the $ax+b$ group. So the Hilbert transformation is invariant under the action of the $ax+b$ group over $L^2(\RR)$. In this sense we say that the Hilbert transformation has the symmetry of the $ax+b$ group. On the other hand this symmetry can characterize the Hilbert transformation. In fact, by using the results in \cite{jo} one first obtains that in $\RR^n,$ and in some other symmetric manifolds as well, linear operators on $L^p$ for some $1\<p<\infty$ commutating with translation or scaling like operations are themselves bounded operator. As second step in \cite{king}
 the author shows that a linear bounded operator on $L^p(\RR)$, $1\<p<\infty$, commuting with the translation and scaling, is of the form
 \begin{equation}\label{}
   \lambda I+\eta {\bf H},
 \end{equation}
 where ${\bf I}$ and ${\bf H}$ are the identity operator and Hilbert transformations respectively, $\lambda, \eta$ are complex numbers.

 \vspace{ 3mm}

\noindent In the present paper by using Gelfand-Naimark's irreducible representation of the ax+b group and Schur theorem we prove the characterization of the Hilbert transformation in the $L^2({\RR})$ space in terms of the $ax+b$ group. We, in particular, do not assume the boundedness property of the operator commuting with the $ax+b$ group. With our approach the boundedness is a consequence of the Gelfand-Naimark's irreducible representation theorem that avoids the big machinery establishes in \cite{jo}.

\vspace{3mm}

\noindent It is more delicate to study the case of  the unit circle $\TT$. Denote by $L^2(\TT)$  the space of square integrable functions on $\TT$ with the inner product $\f1 {2\pi}\int_0^{2\pi}f(e^{i\theta})\ol{g(e^{i\theta})}d\theta $.
 The Hilbert transformation over $L^2(\TT),$ or \emph{circular Hilbert transformation}, is defined as
\begin{equation}\label{cauchy}
  \tilde{\mathbf{H}}f(t)={\rm
  p.v.}\f1{2\pi}\int_0^{2\pi}f(e^{is})\cot(\f{\theta-s}2)ds.
\end{equation}
 where $t=e^{i\theta},\tau=e^{is}$.
A closely related singular integral operator is
\begin{equation}\label{cauchy'}
\mathbf{C}f(t)={\rm
p.v.}\f1{2\pi i}\int_\TT\f{f(\tau)}{\tau -t}d\tau={\rm
p.v.}\f1{2\pi}\int_0^{2\pi}\f{f(e^{is})}{e^{is}-e^{i\theta}}e^{is}ds.
\end{equation}
 In this article we call $\mathbf{C}$ \emph{the singular Cauchy transformation}.
 Denote by $\mathbf{H}_0$ the functional
 $\mathbf{H}_0f=\f1{2\pi}\int_0^{2\pi}f(e^{is})ds$ giving rise to the $0$-th Fourier coefficient of the function to be expanded. It is easy to
 check (the Plemelj formula, see, for instance, \cite{miklin} or \cite{qian1}) that
\begin{equation}\label{H-cauchy}
\mathbf{C}=\frac{i}{2}\tilde{\mathbf{H}}+\frac{1}{2}\mathbf{H}_0.
\end{equation}

 $\tilde{\mathbf{H}}$ has properties analogous with $\mathbf{H}$. For instance, $L^2(\TT)$ is the
 direct sum of the two Hardy spaces on, respectively, the two areas of the complex plane divided by the unit circle.
 $\tilde{\mathbf{H}},$ modulo a constant multiple of $\mathbf{H}_0,$ is a linear combination of the projections over the two Hardy spaces, respectively, namely,
 $$\mathbf{I}=\mathbf{P}^++\mathbf{P}^-,\quad \mathbf{P}^\pm =\frac{1}{2}(\mathbf{I}\pm \mathbf{C}), $$
 and
 \[ \tilde{\mathbf{H}}=(-i)(\mathbf{P}^+-\mathbf{P}^-)+i\mathbf{H}_0, \quad \tilde{\mathbf{H}}^2=-\mathbf{I}+\mathbf{H}_0.\]
 With the Fourier expansion
 \[ f(e^{it})=\sum_{k=-\infty}^\infty c_ke^{ikt}\] there, in fact, hold
 \[ \mathbf{P}^+f(e^{it})=\sum_{k=0}^\infty c_ke^{ikt}, \quad \mathbf{P}^-f(e^{it})=\sum_{k=-\infty}^{-1} c_ke^{ikt},\]
 and
 \[ \tilde{\mathbf{H}}f(e^{it})=-ic_0+\sum_{k\neq0} (-i{\rm sgn}(k))c_ke^{ikt}.\]

\vspace{2mm}

\noindent It seems to be unnatural to study symmetry of the singular Cauchy transformation due to the non-zero curvature of the underlying manifold, viz., the circle. In the present paper we are to deal with symmetry of the circular Hilbert transformation.
 At the first glance, it should be the M\"obius transformation group that gives rise to the characterization of  $\tilde{\mathbf{H}}$.
  There, however, does not seem to exist a Fourier correspondence of the M\"obius transformation.
  On the other hand, the phase translation and scale change generate the Fourier inverse of the actions of the
  $ ax+b$ group  on $L^2(\RR)$.  We were also to obtain the symmetry by the module of the $ax+b$
  group, and, in order to do so, we treat $(a,b)$ and $(c,d)$ as identical if $ax+b\equiv cx+b \ ({\rm mod} 2\pi) \forall x\in \RR.$ Unfortunately,
  the equivalent classes do not form a group.

We construct a family of transformations over $\TT$ whose natural
representation is irreducible over $L^2(\TT)$. Then we obtain the
characterization of  $\tilde{\mathbf{H}}$ in analogy with $\mathbf{H}$.

 We at the end of the paper add some remarks on the role of the M\"{o}bius group in relation to symmetry of the Hilbert transformation on the unit circle.


\section{Induced representations of two groups}

In Section 1 we mentioned that translations and dilations generate  a nontrivial group $G$, the $ax+b$ group,  which is the
group of all affine transformations $x\rt ax+b$ of $\RR$ with $a>0$
and $b\in \RR$. Its underlying manifold is $(0,\?)\times\RR$ and the
group law is defined by
$$
((a,b)(a',b'))(x)=aa'x+b+ab'=(aa',b+ab')(x),\quad  x\in\RR,
$$
which gives
\begin{equation}\label{mult}
\nonumber   (a,b)(a',b')=(aa',b+ab').
\end{equation}
It is  easy to check the relation $(a,b)^{-1}=(\f1a,-\f ba)$. The measure $da/a^2$ is the
left Haar measure and $dadb/a$ is the right Haar measure on this
group \cite{folland}.

There exists a natural unitary representation $\pi$ of $G$, of infinite
dimension, over the Hilbert space $L^2(\RR)$. Denote by
$\mathfrak{U}(L^2(\RR))$ the operator group of the unitary
automorphism of $L^2(\RR)$. Then the group morphism
$\pi:G\rightarrow\mathfrak{U}(L^2(\RR))$ is defined by
\begin{equation}\label{nature1}
 \nonumber   (\pi(a,b)f)(x)=(\f1a)^{\f12}f(\f{x-b}a),\ x\in\RR.
\end{equation}
$\pi(a,b)$ is also written as $\pi_{ab}$ in this article.
\vskip2mm
It is obvious that  both the Hardy spaces on respectively the upper and lower half planes are the invariant subspaces of
 the natural representation of the $ax+b$ group. So it is reducible. It is not easy to get the irreducible
 representation of the $ax+b$ group since it is both noncommutative and noncompact. In 1948, Gelfand and Naimark  \cite{gefand} first proved the following Theorem.

\bT[Gelfand-Naimark]\label{gelfand} The $ ax+b$ group has only two
nontrivial irreducible representation,
 $\check{\pi}^+(a,b) :(L^2(0,+\infty)\to L^2(0,+\infty))$ and $ \check{\pi}^-(a,b): (L^2(-\infty,0) \to L^2(-\infty,0))$
 as
 \begin{eqnarray}
 \label{g1} &[\check{\pi}^+(a,b)f](x)=a^{\frac{1}{2}}e^{2\pi ibx}f(ax)  ,(x>0) \\
 \label{g2} &[\check{\pi}^-(a,b)f](x)=a^{\frac{1}{2}}e^{2\pi ibx}f(ax)  ,(x<0)
 \end{eqnarray}
 \eT

\noindent We define a representation of the $ax+b$ group over $L^2(\RR)$ as
\begin{equation}\label{rrep}
 ( \check{\pi}(a,b)f)(x)=a^{\frac{1}{2}}e^{2\pi ibx}f(ax),\ a\in\RR^+,b\in\RR, f\in L^(\RR).
\end{equation}
$\check{\pi}^\pm(a,b)$ and $ \check{\pi}(a,b)$ are also written as $\check{\pi}^\pm_{ab}$ and  $ \check{\pi}_{ab}$, respectively.
\vskip2mm
  It is obvious that $L^2(\RR)$ can be decomposed into the orthogonal direct
 sum
 of $L^2(-\infty,0)$ and $L^2(0,+\infty)$, i.e,
$L^2(\RR)=L^2(-\infty,0)\oplus L^2(0,+\infty)$. By Theorem
\ref{gelfand},  the  representation $\check{\pi}$ is just  the sum of two
irreducible representations $\check{\pi}^\pm$ over $L^2(-\infty,0)$
and $L^2(0,+\infty)$ respectively.

 Denote by $\mathcal{F}$ the
Fourier transformation. We also denote $f^{\wedge}=\mathcal{F}f.$  For $(a,b) $ in the $ax+b$ group,
define $\mathcal{F}({\pi})$ by
\begin{equation}\label{}
   ((\mathcal{F}({\pi})({a,b}))f)(x)
  =(\mathcal{F}^{-1}({\pi}_{ab}f^{\wedge}))(x),\ f\in L^2{(\RR)}.
\end{equation}
Then we have
\begin{eqnarray}
  \nonumber ((\mathcal{F}({\pi})({a,b}))f)(x)   &=& \f1{\sqrt{2\pi}}\int_{-\?}^\?a^{-\frac{1}{2}}f^{\wedge}(\f y a-\f b
   a)e^{iyx}dy\\
   \nonumber &=&a^{\frac{1}{2}}\f1{\sqrt{2\pi}}\int_{-\?}^\?f^{\wedge}(\f {y -
   b}
   a)e^{i\f{y-b}a(ax)}e^{ibx}d\f{y-b}a\\
&=&a^{\frac{1}{2}}e^{ibx}f(ax)=(\check{\pi}_{a,b}f)(x).
\end{eqnarray}
 Thus the Fourier transformation is just the isomorphism between the two representation $\pi$
with $\check{\pi}$.
Then $\pi$ has two sub-representations that are equivalent to $\check{\pi}^\pm$ by the Fourier correspondence, respectively. Let $H^\pm$ be the the images of
 $L^2(0,+\infty)$ and $L^2(-\infty,0)$ under the Fourier transformation, respectively. Let $\pi^\pm$ be the restriction of $\pi$ over $H^\pm$, respectively.
 \bT\label{pi-rep} For $\pi=\pi_{ab},$ we have $\pi=\pi^+\bigoplus\pi^-$, where $\pi^\pm$ are defined by
\begin{eqnarray}
\nonumber{\pi}^+f&=&(\check{\pi}^+(a,b)f^{\vee})^{\wedge},\ f\in H^+;\\
  \nonumber  {\pi}^-f&=&(\check{\pi}^-(a,b)f^{\vee})^{\wedge}, \ f\in H^-.
    \end{eqnarray}
Then $\pi^\pm$ are irreducible representations of the $ax+b$ group over $H^\pm$, respectively. \eT

\noindent To obtain our main results we introduce a Schur's lemma in the version of the infinite dimension. Let
$\mathcal{H}$ be a complex Hilbert space, and $\mathcal{A}$  a
family of transformations acting on $\mathcal{H}$ satisfying $\mathcal{H}=\bigcup_{T\in \mathcal{A}}T(\mathcal{H})$. We say that
$\mathcal{A}$ acts irreducibly on $\mathcal{H}$ if there does not
exist a decomposition of $\mathcal{H}$ such that
$\mathcal{H}_1\oplus\mathcal{H}_2,$ where both $\mathcal{H}_1$ and
$\mathcal{H}_2$ are invariant subspaces of $\mathcal{A}$.

\bT\cite{kn} {\bf(Dixmier's Lemma)}\label{schur} Let $\mathcal{H}$ be a  complex separable Hilbert space. If a
family of transformations $\mathcal{A}$ acts irreducibly on
$\mathcal{H}$, then any linear transformation $\mathbf{T}$ which commutes
with every $T\in \mathcal{A}$ is of the form $\mathbf{T}=\lm\mathcal{I}$ where $\lm$ be
a complex number and $\mathcal{I}$ is the identity transformation.\eT

Let $\s$ be a representation of a group $G$ over the Hilbert space
$\mathcal{H}$. It is obvious that $\s$ is  irreducible
  if and only if $\{\s(g):g\in G\}$ acts irreducibly on
$\mathcal{H}$. Then by the above theorem, the following theorem is evident.

\bT {\bf(Schur's Lemma)}\label{t1} Suppose that $\s$  is a irreducible
representation of a group $G$ over a complex separable Hilbert space $\mathcal{H}$. Then any linear transformation $\mathbf{T}$ which commutes
with $\s$ is of the form  $\mathbf{T}=\lm\mathcal{I}$ where $\lm$ is
a complex number and $\mathcal{I}$ is the identity transformation.\eT
We note that $\mathbf{T}$ is not required to be a bounded operator in either of the above theorems.

\noindent We will be based on the above two theorems.
 \section{Characterization of the Hilbert transformation on the Line}
We can now state the symmetry properties of the Hilbert transformation $\mathbf{H}$. It is obvious that the operator group $\pi(G)$ is  generated by   $T_a$ and $\tau_b$. Then by Lemma \ref{commute} we get that the Hilbert transformation $\mathbf{H}$ is invariant under the actions of $\pi(G)$ over $L^2(\RR)$. That is the following theorem.

 \bT\label{pi-com}$\pi$ commutes with the Hilbert transform $\mathbf{H}$, i.e.,
 \begin{equation}\label{commu}
   \pi_{ab}\mathbf{H}(f)=\mathbf{H}\pi_{ab} (f),\quad (a,b)\in \mathbf{ax+b},f\in L^2(\RR).
 \end{equation}

 \eT

Then by Theorem \ref{pi-rep} and \ref{t1}, the restrictions of $\pi$ over $H^\pm$ are the scalar operators. Let us now further identify the spaces $H^\pm$. Denote by $H^2_{\pm}$ the Hardy spaces on the upper and lower half planes respectively.
\bT \label{hardy}$H^+=H^2_{+}$ and $H^-=H^2_{-}$. \eT
\bpp
It is trivial to check that both $H^2_{\pm}$ are the closed invariant subspaces of $\pi(G)$. By Theorem \ref{gelfand} there exist exactly two nontrivial irreducible representations of the $ ax+b$ group. Then $H^2_{\pm}$ must be $H^\pm.$

 Given $0\ne f\in H^+$. By the definition of $H^+$, its Fourier transform $f^{\wedge}\in L^2(0,+\infty)$. We define
 $$F(z)=\frac{1}{2\pi}\int_{-\infty}^{+\infty}f^{\wedge}(t)e^{itz}\,dt=\frac{1}{2\pi}
   \int_{0}^{+\infty}f^{\wedge}(t)e^{itx}e^{-ty}\,dt,   \ z=x+iy.$$
Obviously, for $y > 0,$ $F(z)$ is analytic. By the theory of the Hardy space on the upper half plane we obtain that $f(x)$ is the non-tangential boundary
limit of $F(z)$ when $z$ tends to $x$ from the above. That is,
$f(x)$ is the boundary value of $F(z)$, i.e, $f$ belongs to $H^2_{+}$. So $H^+=H^2_{+}$, and then  $H^-=H^2_{-}$.\epp

We note that from the Fourier multiplier representation of the Hilbert transformation we have
\begin{equation}\label{hs}
  ({\mathbf{H}f})^\wedge(\xi)=-i{\rm sgn}(\xi)f^\wedge(\xi),\ \xi\in\RR.
  \end{equation}

The following theorem is now obvious.

\bT\label{char} $H^+$ and $H^-$ are respectively the eigen-subspaces associated to eigenvalues $- i$ and $i$ of the Hilbert transformation. That is
$$\mathbf{H}|_{H^+}=-i\mathbf{I}|_{H^+},H|_{H^-}=i\mathbf{I}|_{H^-}.$$\eT
 Although the Hardy space theory and the related Fourier multiplier theory imply Theorem \ref{char}, it is difficult for them to discuss the converse of the theorem. The converse of the theorem addresses the symmetry of the Hilbert transformation.
 \bT
 Suppose that $\mathbf{T}$ is a linear operator from $L^2(\RR)$ to itself, and $\mathbf{T}$
   commutes with the natural representation of the group $ax+b.$ Then there exist two complex numbers
 $\lambda,\eta$ such that
 \begin{equation}\label{inverse}
   \mathbf{T}=\lambda \mathbf{I}+\eta \mathbf{H}.
 \end{equation}
Moreover, if $\mathbf{T}$ is an anti-symmetric, norm-preserving and real operator, it must be either $\mathbf{H}$ or $-\mathbf{H}$. \eT
  \bpp
  By Theorem \ref{t1} both the restrictions of  $\mathbf{T}$ over $H^\pm$ are the scalar operators since $\mathbf{T}$ commutes with $\pi(G)$.
Assume that $\mathbf{T}=k_1 \mathbf{I}|_{H^+}$ over $H^+$ and $\mathbf{T}=k_2 \mathbf{I}|_{H^-}$ over $H^-$, respectively. For $f\in L^2(\RR)$  there exist $f_1\in H^+$ and $f_2\in H^-$ such that $f=f_1+f_2$. Then we have
\begin{eqnarray}\label{proof of inverse}
    \nonumber\mathbf{T}f &=&  \mathbf{T}f_1+ \mathbf{T}f_2=k_1f_1+k_2f_2\\
   \nonumber &=& \f{k_1+k_2}2(f_1+f_2)+\f{k_1-k_2}{2i} (if_1-if_2)\\
   &=&  \f{k_2+k_1}2\mathbf{I}f+\f{k_2-k_1}{2i} \mathbf{H}f.
\end{eqnarray}
Let $\lambda=\frac{k_2+k_1}{2}$ and $\eta=\frac{k_2-k_1}{2i}$. (\ref{proof of inverse})
completes the proof of (\ref{inverse}).

If $\mathbf{T}$ is a real operator, its non-real eigenvalues must appear in pairs. This implies $k_1=\overline{k_2}$.
If $\mathbf{T}$ is also anti-symmetric, its eigenvalues must be pure complex numbers, i.e., $k_1=-k_2$. Finally, if  $\mathbf{T}$ preserves the norm, we obtain that $|k_1|=|k_2|=1$. Therefore, $\lambda=0$ and $\eta=1$ or $-1$.
\epp
\vskip2mm

\brs  For bounded operators $T$ Theorem \ref{char} is a  known result in the literature of Hilbert transformation. The above (\ref{proof of inverse}) is precisely the relation (4.82) in \cite{king}. Moreover, linear self-maps of $L^2(R)$ that commute with translations (even just with one translation) are automatically continuous \cite{jo}. Since translation operators over $L^2(\RR)$ do not have critical eigenvalue, by Corollary 3.5 of \cite{jo}, linear self-maps of $L^2(R)$ that commute with translations must be continuous. The relation (4.82) in \cite{king} is proved by the multiplier theory from \cite{stein}, where the assumption of boundedness of $T$ is essential. In the present paper we derive Theorem \ref{char} without invoking the results from \cite{jo,king,stein}.
 We use the Gelfend-Naimark¡¯s representation of $ax+b$ to deal with the issue. The new methodology would enable us to study  Hilbert transformations on other types of manifolds with symmetry properties similar to translations and dilations. As example, in the next section we study the Hilbert transformation on the unit circle in an analogous way.
\ers

 \section{Characterization of the Hilbert transformation on the Circle}

We did not find a group acting on $\TT$ suitable to study the Hilbert transformation over $L^2({\TT})$. Fortunately, there exists a semigroup on the unit circle that plays a similar role as the
$ax+b$ group on the real line. The semigroup on the unit circle makes it possible to
carry on a similar but not exactly the same procedure as above for the line.

Denote by $\mathfrak{N}$ the set of pairs $(n,\beta),n\in\ZZ^+, \beta\in \RR$. We will also use the notation $\a\theta+\b$, $\a\in\QQ^+,\b\in\RR$. The latter turns out to be a subgroup of the  $ax+b$ group.  Be
similar to the case of the $ax+b$ group, we would hope that  $\a\theta+\b$
can be considered as an affine transformation on $\RR$ such that
$(n,\beta)(\theta)=n\theta+\beta\ \mod(2\pi)$. Let $(n,\beta)\in \mathfrak{N}$. Define its action on $\TT$ by $(n,\beta)t=e^{i(n\theta+\b)}$, where $t=e^{i\theta}\in \TT$.  Unfortunately it
does not work for the non-integer $\alpha$'s by this machine.
That is, $\a\theta+\beta$ is not
congruent to $\a(\theta+2\pi)+\beta$ modulo $2\pi$ if $\a$ is not an integer.  It means that
it is not well-defined for the action of $\a\theta+\b$  on $\TT$.
\vskip2mm

In general the element of the  ${\bf \a\th+\b}$ group is not the automorphism of $\TT$. There, however, still exists  a natural action of the \bf $\a\th+\b$ \rm group over the space
$L^2(\TT)$.

Define by $\pi$ the mapping from the $\a\th+\b$  group to the
the set of bounded endomorphisms of $L_2(\TT)$  as follows: for
$n\in \ZZ^+$ and $f\in L^2(\TT),$
\begin{equation}\label{integer}
    (\pi(n,\b))f(t)=(\f1{n})^{\f12} f(t^ne^{i\b}),\ t\in\TT;
\end{equation}
and, for a general positive rational number $\a$, $\a=q/p, \ p, \ q\in
\ZZ^+,\ p>1$ and $(p,q)=1$, let
\begin{equation}\label{nointeger}
  (\pi(\f q p,\b)f)(t)= (\f p q)^{\f12}\f1p(f(e^{i(\f q p\th+\b)})+f(e^{i(\f q p\th+\b)}\om_p^1)+\cdots+f(e^{i(\f q p\th+\b)}\om_p^{p-1})),
   \end{equation}
where $\t=e^{i\th}\in\TT$ and $\omega_p=e^{i\f{2\pi}p}$.

Notice that
$e^{i\f{q2k\pi}p}$ is one of the $p$-th roots of the unity. Then the $p$-tuple\\
 $(e^{i\f{q2k\pi}p},e^{i\f{q2k\pi}p}w_p^1, \cdots,e^{i\f{q2k\pi}p}w_p^{ {p}-1})$ is just a
rearrangement of $(1,\om^1_{ {p} },\cdots,\om_{p}^{p-1})$. Thus (\ref{nointeger}) is independent of the choice of $\th$.
We note that (\ref{integer}) and (\ref{nointeger}) do not give the representation of the $\a\th+\b$ group. But we still can use it to derive the characterization of the Hilbert transformation over $L^2(\TT)$.

In the rest of the section we denoted by  $\pi_{\a\b}$ the image of
$(\a,\b)$ under $\pi$. We say that $\pi$ commutes with a linear operator $\mathbf{T}$ if $\pi_{\a\b}\mathbf{T}=\mathbf{T}\pi_{\a\b}$ for all $(\a,\b)$ in the $\a\th+\b$ group.

\bT\label{main} $\pi$ commutes with the Hilbert transformation on the
unit circle.\eT
 \bpp  Let $p,q\in\ZZ$ and $(p,q)=1$. First  we check that $\mathbf{H}_0$ is invariant under $\pi$. In fact, for $f\in L^2(\TT)$ we have
 \begin{eqnarray}\label{}
   \nonumber (\mathbf{H}_0 (\pi(\f q p,\b)f))(t)&=&\frac{1}{2\pi}\int_0^{2\pi}(f(e^{i(\f q p s+\b)})+\cdots+f(e^{i(\f q p s+\b)}\om_p^{p-1}))ds\\
    \nonumber &=&(\f p q)^{\f12}\frac{1}{2\pi}\int_0^{2\pi}f(e^{is})ds=(\pi(\f q p,\b)(\mathbf{H}_0f))(t).
\end{eqnarray}
 Then by (\ref{H-cauchy}) it is sufficient to prove that $\pi$ commutes with the singular Cauchy transformation $\mathbf{C}$.
 Notice that
$\pi_{\f q p \b}$ is the composition of $\pi_{q0}$ with $\pi_{\f 1p
\b}$, i.e., $\pi_{\f q p \b}=\pi_{q0}\pi_{\f 1p \b}$. We only prove
the cases of $\a=n$ or $\f1n$ for positive integers $n$.

Let $\a=n$ be a positive integer and $\b\in\RR$. It is obvious that
we may assume $n>1$. Then for $f\in L^2(\TT)$ with the H\"older type continuity we have
\begin{eqnarray}\label{htcom}
    (\mathbf{C}\pi_{n\b})f(s)&=&\f1{2\pi i}{\rm p.v.}\int_\TT\f{(\f1{n})^{\f12}f((n,\b)t)}{t-s}dt\\
   \nonumber &=&\f1{2\pi i}{\rm p.v.}\int_\TT \f{(\f1{n})^{\f12}f(t^ne^{i\b})}{t-s}dt\\
   \nonumber &=&(\f1{n})^{\f12}\f1{2\pi
    }{\rm p.v.}\int_0^{2\pi}\f{f(e^{i(n\th+\b)})}{e^{i\th}-s}e^{i\th}{d\th}
\end{eqnarray}
Let $\phi=n\th$. Then by (\ref{htcom}) we get
\begin{eqnarray}\label{htcom1}
    (\mathbf{C}\pi_{n\b})f(s)&=&(\f1{n})^{\f12}\f1{2\pi
    }{\rm p.v.}\int_0^{2n\pi}\f{f(e^{i(\phi+\b)})}{e^{i\f 1n\phi}-s}e^{i\f1n \phi}\f{d\phi}{n }\\
  \nonumber &=&(\f1{n})^{\f12}\f1{2\pi
    }{\rm p.v.}\int_0^{2n\pi}\f{f(e^{i(\phi+\b)})}{e^{i\phi}-s^n}e^{i\f1n
    \phi}(e^{i\f{n-1}n\phi}+\cdots+s^{n-1})\f{d\phi}{n }\\
  \nonumber  &=&\f1n(\f1{n})^{\f12}\f1{2\pi
    }({\rm p.v.}\int_0^{2n\pi}\f{f(e^{i(\phi+\b)})e^{i\phi}}{e^{i\phi}-s^n}{d\phi}+\cdots\\
\nonumber&\
&+{\rm p.v.}\int_0^{2n\pi}\f{f(e^{i(\phi+\b)})s^{n-1}}{e^{i\phi}-s^n}{d\phi})
\end{eqnarray}
We divide the right hand side of (\ref{htcom1}) into $n$ parts. Then
\begin{eqnarray}\label{1}
\f1n(\f1{n})^{\f12}\f1{2\pi
    }{\rm p.v.}\int_0^{2n\pi}\f{f(e^{i(\phi+\b)})e^{i\phi}}{e^{i\phi}-s^n}{d\phi}&=&
(\f1{n})^{\f12}\f1{2\pi
    }{\rm p.v.}\int_0^{2\pi}\f{f(e^{i(\phi+\b)})e^{i\phi}}{e^{i\phi}-s^n}{d\th}\\
 \nonumber&=&(\f1{n})^{\f12}\f1{2\pi i
    }{\rm p.v.}\int_0^{2\pi}\f{f(e^{i\psi})de^{i\psi}}{e^{i\psi}-s^ne^{i\b}}\\
\nonumber&=& (\f1n)^{\f12}
 \f1{2\pi
i}{\rm p.v.}\int_\TT\f{f(\tau)d\tau}{\tau-s^ne^{i\b}}\\
\nonumber&=&\pi_{n\b}(\mathbf{C}f)(s).
\end{eqnarray}
Similarly, for $1<k\<n-1$ we have
\begin{eqnarray}\label{2}
   {\rm p.v.} \int_0^{2n\pi}\f{f(e^{i(\phi+\b)})e^{i\f {n-1-k}n\th}s^{k}}{e^{i\phi}-e^{i\b}}{d\th}
   &=&
     \sum_{l=1}^n{\rm p.v.}\int_{2(l-1)}^{2l\pi}\f{f(e^{i(\phi+\b)})e^{i\f {n-1-k}n\th}s^{k}}{e^{i\phi}-s^n}{d\th} \\
 \nonumber  &=& {\rm p.v.}\int_{0}^{2\pi}\f{f(e^{i(\phi+\b)})s^{k}}{e^{i\phi}-e^{i\b}}(\sum_{l=1}^ne^{i\f
 {n-1-k}n(\th+2l\pi)}){d\th}\\
  \nonumber&=& {\rm p.v.}\int_{0}^{2\pi}\f{f(e^{i(\phi+\b)})s^{k}}{e^{i\phi}-e^{i\b}}e^{i\f {n-1-k}n\th}(\sum_{l=1}^ne^{i\f
  {n-1-k}n(2l\pi)}){d\th}.
\end{eqnarray}
Notice that $n>k>1$ and $e^{i\f
  {n-1-k}n(2l\pi)},l=1,2,\cdots,n,$ are just all of the $n$-roots of
 the  unity, which means that
\begin{equation}\label{3}
\sum_{l=1}^ne^{i\f
  {n-1-k}n(2l\pi)}=0.
\end{equation}
Now by (\ref{htcom1}),(\ref{1}),(\ref{2}) and (\ref{3}), the
proposition is proved for $\a\in\ZZ^+$.

For $\a=\f1n$ we have
\begin{eqnarray}\label{htcomf}
    (\mathbf{C}\pi_{\a\b})f(s)&=&\f1{2\pi i}{\rm p.v.}\int_\TT\f{(\pi_{\a\b}f)(t)}{t-s}dt\\
   \nonumber &=&\f1{2\pi i}{\rm p.v.}\int_0^{2\pi} {n}^{\f12}\f1n(f(e^{i\f1n\th}e^{i\b})+\cdots\\
 \nonumber  &\ &+f(e^{i\f1n\th}\om_n^{n-1}e^{i\b}))\f1{e^{i\th}-s}de^{i\th}.
\end{eqnarray}
For $1\<k\<n$, let $\phi=\f1n\th$ in the $k$-th terms of the  sum
of the above integrand. Then we have
\begin{eqnarray}
\label{ht1}&\ &\f1{2\pi i}{\rm p.v.}\int_0^{2\pi }
{n}^{\f12}\f1n(f(e^{i\f1n\th}\om_n^{k-1}e^{i\b})
\f{de^{i\th}}{e^{i\th}-s}\\
\nonumber  &=&\f1{2\pi }
{n}^{\f12}{\rm p.v.}\int_0^{\f2n\pi}\f{f(e^{i\phi}\om_n^{k-1}e^{i\b})e^{in\phi}}{e^{in\phi}-s}d\phi \\
  \nonumber   &=&\f1{2\pi }
{n}^{\f12}{\rm p.v.}\int_0^{\f2n\pi}\f{f(e^{i\phi}\om_n^{k-1}e^{i\b})e^{i(n-1)\phi}e^{i\phi}}
{(e^{i\phi}-e^{i\f1n\psi})(e^{i\phi}-e^{i\f1n\psi}w_n)\cdots(e^{i\phi}-e^{i\f1n\psi}w_n^{n-1})}d\phi,
\end{eqnarray}
where $s=e^{i\psi}\in\TT$ and $w_n=\f1{\om_n}=e^{-i\f{2\pi}n}$. Notice that $w_n$ is
still one of the $n$-th roots of the unity. Then we obtain that $(1-w_n)\cdots(1-w_n^{n-1})=\lim_{z\rightarrow1}\f{z^n-1}{z-1}=n$ and
\begin{eqnarray}\label{frac}
&\ &\f{e^{i(n-1)\phi}}
{(e^{i\phi}-e^{i\f1n\psi})(e^{i\phi}-e^{i\f1n\psi}w_n)\cdots(e^{i\phi}-e^{i\f1n\psi}w_n^{n-1})}\\
 \nonumber &=&\f{1}{(1-w_n)\cdots(1-w_n^{n-1})}(\f1{e^{i\phi}-e^{i\f1
n\psi}}+\cdots+\f1{e^{i\phi}-e^{i\f1 n\psi}w_n^{n-1}})\\
 \nonumber &=&\f1n(\f1{e^{i\phi}-e^{i\f1
n\psi}}+\cdots+\f1{e^{i\phi}-e^{i\f1 n\psi}w_n^{n-1}}).
\end{eqnarray}
Denote by
\begin{equation}\label{aij}
  A_{kj}=\f1{2\pi }
\f{{n}^{\f12}}n{\rm
p.v.}\int_0^{\f2n\pi}\f{f(e^{i\phi}\om_n^{k-1}e^{i\b})e^{i\phi}}{e^{i\phi}-e^{i\f1
n\psi}w_n^{j-1}}d\phi,\ j,k=1,\cdots,n-1.
\end{equation}
From (\ref{htcomf})-(\ref{aij}) we have
\begin{equation}\label{caij}
(\mathbf{C}\pi_{\a\b})f(s)=\sum_{k=1}^n\sum_{j=1}^nA_{kj}=\sum_{m=0}^{n-1}(\sum_{1\<k,j\<n,j-k\equiv
m({\rm mod}\ n)}A_{k,j}).
\end{equation}
It is obvious that $\sum_{1\<k,j\<n,j-k\equiv 0({\rm mod}\
n)}A_{k,j}= \sum_{j=1}^nA_{j,j}$. But by (\ref{aij}) we obtain
\begin{eqnarray}
\label{11} \sum_{j=1}^nA_{j,j}
&=&\f1{2\pi } \f{{n}^{\f12}}n{\rm
p.v.}(\int_0^{\f{2\pi}n}\f{f(e^{i\phi}e^{i\b})e^{i\phi}}{e^{i\phi}-e^{i\f1n\psi}}d\phi+
\int_0^{\f{2\pi}n}\f{f(e^{i\phi}\om_ne^{i\b})e^{i\phi}}{e^{i\phi}-e^{i\f1n\psi}w_n}d\phi\\
\nonumber &\ &+\cdots
\int_0^{\f{2\pi}n}\f{f(e^{i\phi}\om^{n-1}_ne^{i\b})e^{i\phi}}{e^{i\phi}-e^{i\f1n\psi}w^{n-1}_n}d\phi )\\
 \nonumber  &=&\f1{2\pi i} \f{{n}^{\f12}}n{\rm
p.v.}(\int_0^{\f{2\pi}n}\f{f(e^{i\phi}e^{i\b})}{e^{i\phi}-e^{i\f1n\psi}}de^{i\phi}+
\int_0^{\f{2\pi}n}\f{f(e^{i\phi}e^{i\f{2\pi}n}e^{i\b})e^{i\f{2\pi}n}}{e^{i\phi}e^{i\f{2\pi}n}-e^{i\f1n\psi}}
de^{i\phi}\\
\nonumber &\ &+\cdots \int_0^{\f{2\pi
}n}\f{f(e^{i\phi}e^{i\f{2(n-1)\pi}n}e^{i\b})e^{i\f{2(n-1)\pi}n}}{e^{i\phi}e^{i\f{2(n-1)\pi}n}
-e^{i\f1n\psi}}de^{i\phi}) \\
\nonumber &=&\f1{2\pi i} \f{{n}^{\f12}}n{\rm
p.v.}(\int_0^{\f{2\pi i}n}+\cdots+\int_{\f{2(n-1)\pi}n}^{2\pi})\f{f(e^{i\phi}e^{i\b})de^{i\phi}}{e^{i\phi}-e^{i\f1n\psi}}\\
\nonumber &=&\f1{2\pi i} \f{{n}^{\f12}}n{\rm
p.v.}\int_0^{2\pi}\f{f(e^{i\phi}e^{i\b})de^{i\phi}}{e^{i\phi}-e^{i\f1n\psi}}
=\f1{2\pi i} \f{{n}^{\f12}}n{\rm
p.v.}\int_0^{2\pi}\f{f(e^{i\phi})de^{i\phi}}{e^{i\phi}-e^{i\f1n\psi}e^{i\b}}
\end{eqnarray}
Similarly, for $m\geq1$, recalling that
$\om^k_n,\om_n^{k+1},\cdots,\om_n^{k+n-1}$ still run through the
$n$-th roots of the unity, we get
\begin{eqnarray}\label{caijj}
\label{22}   &\ &\sum_{1\<k,j\<n,j-k\equiv m({\rm mod}\
n)}A_{k,j}\\
\nonumber&=&\f1{2\pi } \f{{n}^{\f12}}n{\rm
p.v.}(\int_0^{\f{2\pi}n}\f{f(e^{i\phi}\om^m_ne^{i\b})}{e^{i\phi}-e^{i\f1n\psi}}de^{i\phi}+
\int_0^{\f{2\pi}n}\f{f(e^{i\phi}\om^{m+1}_ne^{i\b})}{e^{i\phi}-e^{i\f1n\psi}w_n}de^{i\phi}\\
\nonumber&\
&+\cdots+\int_0^{\f{2\pi}n}\f{f(e^{i\phi}\om^{n-1}_ne^{i\b})}{e^{i\phi}-e^{i\f1n\psi}w_n^{n-m-1}}de^{i\phi}
+\int_0^{\f{2\pi i}n}\f{f(e^{i\phi}e^{i\b})}{e^{i\phi}-e^{i\f1n\psi}w_n^{n-m}}de^{i\phi}\\
\nonumber&\ &+\cdots+\int_0^{\f{2\pi}n}\f{f(e^{i\phi}\om^{m-1}_ne^{i\b})}{e^{i\phi}-e^{i\f1n\psi}w^{n-1}_n}de^{i\phi} )\\
\nonumber&=&\f1{2\pi i} \f{{n}^{\f12}}n{\rm
p.v.}(\int_0^{\f{2\pi}n}\f{f(e^{i\phi}\om^m_ne^{i\b})}{e^{i\phi}-e^{i\f1n\psi}}de^{i\phi}+
\int_0^{\f{2\pi}n}\f{f(e^{i\phi}\om^{m+1}_ne^{i\b})}{e^{i\phi}-e^{i\f1n\psi}w_n}de^{i\phi}\\
\nonumber&\
&+\cdots+\int_0^{\f{2\pi}n}\f{f(e^{i\phi}\om^{n-1}_ne^{i\b})}{e^{i\phi}-e^{i\f1n\psi}w_n^{n-m-1}}de^{i\phi}
+\int_0^{\f{2\pi}n}\f{f(e^{i\phi}\om_n^n e^{i\b})}{e^{i\phi}-e^{i\f1n\psi}w_n^{n-m}}de^{i\phi}\\
\nonumber&\ &+\cdots+\int_0^{\f{2\pi}n}\f{f(e^{i\phi}\om^{n+m-1}_ne^{i\b})}{e^{i\phi}-e^{i\f1n\psi}w^{n-1}_n}de^{i\phi} )\\
\nonumber&=&\f1{2\pi i} \f{{n}^{\f12}}n{\rm
p.v.}\int_\TT\f{f(\tau)}{\tau-e^{i\f1n\psi}\om_n^me^{i\b}}d\tau,
\quad m=1,2,\cdots,n-1.
\end{eqnarray}
 From (\ref{caij}) and  (\ref{caijj}),  we have
\begin{equation}\label{end}
    (\mathbf{C}\pi_{\f1n,\b})f(s)=n^{\f12}\f1n(\sum_{m=0}^{n-1}\f1{2\pi
    i}\int_\TT\f{f(\tau)}{e^{i\f1n\psi}\om_n^{m}e^{i\b}-\tau}d\tau)=\pi_{\f1n,\b}(\mathbf{C}f){s}.
\end{equation}
(\ref{end}) is valid for $f\in L^2(\TT)$ because  $\mathbf{C}$ is a
bounded operator over $L^2(\TT) $ and the class of the functions of the H\"older continuity is dense
in $L^2(\TT) $. This completes the proof of the theorem.

 \epp

\section{Decomposition of the concerned operators and spaces on the citcle}

In the rest of the paper we discuss the irreducibility of
$\pi$. Then we  characterize the circular Hilbert transformation $\tilde{ \mathbf{H}}.$

In the previous section we point out that   $G=\{(\a,\b),\a\
\text{or}\ \f1\a\in \ZZ^+,;b\in \RR\}$ is not a group. But  we
can still prove that the family of $\pi_{(\a,\b)},\ (\a,\b)\in G,$ acts
irreducibly on some subspaces of $L^2(\TT)$.

Given a family of functions $\mathfrak{M}$ in  $L^2(\TT)$. Denote by
\begin{equation}\label{zero}
  Z(\mathfrak{M})=\bigcap_{f\in\mathfrak{M}}\{n\in \ZZ,f_n=0\},
\end{equation}
where $f_n$ are the Fourier coefficients of $f$ such that
$f(t)=\sum\limits_{n=-\?}^\?f_nt^n$. By (\ref{zero}),
$Z(\mathfrak{M})=\emptyset$ means that for every $n\in\ZZ$ there
exists at least one $f^{(n)}\in \mathfrak{M}$ such that
$f^{(n)}_n\neq0$. \bT\label{con0} Assume that $ Z(\mathfrak{M})$ is
empty. Suppose that $\phi$ belongs to $ L^2(\TT)$ and satisfies
\begin{equation}\label{con}
  f\ast \phi\equiv0,\quad \forall f\in \mathfrak{M}.
\end{equation}
Then $\phi\equiv0$.\eT
\bpp
Assume that
\begin{equation}\label{}
  \nonumber \phi(t)=\sum_{n=-\?}^\?\phi_n t^n=\sum_{n=-\?}^\?\phi_n e^{in\theta},\quad t=e^{i\theta};
\end{equation}
and
\begin{equation}\label{}
  \nonumber f(t)=\sum_{n=-\?}^\?f_n t^n=\sum_{n=-\?}^\?f_n e^{in\theta},\quad t=e^{i\theta}
\end{equation}
for $f\in\mathfrak{M}$.

Then we obtain
\begin{equation}\label{con1}
  (f\ast \phi)(e^{i\theta})=\mathcal{F}^{-1}(\hat{f}\hat{\phi})=\sum_{n=-\?}^\?f_n\phi_n e^{in\theta},
\end{equation}
which gives $f_n\phi_n=0,\forall f\in\mathfrak{M}$, especially
$f^{(n)}_n\phi_n=0$. Thus $\phi_n=0,\forall n\in \ZZ$, which means
that $\phi\equiv0$. \epp

Denote by $H^+(\TT)$ the Hardy space on the unit disc,  $H^-(\TT)$
the Hardy space on the complement of the unit disc in the whole complex
plane.
  Denote by $H^0$  the subspace of constant functions, and $\tilde{H}^+(\TT)$ its orthogonal
  complement in $H^+(\TT)$. Then we obtain that   $L^2(\TT)=\tilde{H}^+\bigoplus H^0\bigoplus H^-$. It is obvious that
  $H^0$ is invariant to $\pi$. It is easy to check
  that there do not exist $f\in\tilde{H}^+(\TT)$ and  $\pi_{(\a,\b)}$ such that $\pi_{(\a,\b)}f$ is
   a constant function. Thus both $H_0$ and $\tilde{H}^+(\TT)$ are invariant spaces
   of $\pi_{(\a,\b)},\ (\a,\b)\in G$.
\vskip2mm
\noindent The following theorem plays the same role on the circle as the Gelfend-Naimark's representation for the \bf ax+b \rm group on the real axis.
\bT \label{dim} The family of $\pi_{(\a,\b)},\ (\a,\b)\in G,$ acts
irreducibly over  $\tilde{H}^+(\TT)$ and $H^-(\TT)$ respectively.\eT
\bpp We only prove the first part. The proof of the other part is similar.

If the family $\pi_{(\a,\b)}$ were not
irreducible over  $\tilde{H}^+(\TT),$ there would exist two  proper subspaces
$H_1,H_2\subset \tilde{H}^+(\TT)$ such
 that $\tilde{H}^+(\TT)=H_1\bigoplus H_2$, and both $H_1$ and $H_2$ are the
  invariant spaces of $\{\pi_{(\a,\b)},(\a,\b)\in G\}$.

\bf Claim: \rm $Z(H_1)\cap \ZZ^+$ can not be empty.

Otherwise, assume that $Z(H_1)\cap \ZZ^+=\emptyset$. Since
$Z(H^-(\TT))=\ZZ^+\cup \{0\}$, $Z(H^0)=\ZZ\setminus \{0\}$, we
obtain that $Z(H^-(\TT))\cap Z(H^0)=\ZZ^+$ and
\begin{equation}\label{con2}
  Z(H_1\cup H^0\cup H^-(\TT))=Z(H_1)\cap Z(H^0)\cap Z(H^-(\TT))=\emptyset.
\end{equation}
Let $h\in H_2$.  Since all $H^0$, $H_1$ and $H_2$ are invariant
spaces of the transforms $\pi_{(1,\theta)}$ for $\theta\in \RR$ we
have
\begin{equation}\label{}
 \nonumber <\pi_{(1,\theta)}f,h>=0,\forall f\in H_1\cup H^0\cup H^-(\TT),\theta\in
 \RR.
\end{equation}
Recalling that $\pi_{(1,\theta)}f(e^{is})=f(e^{i(s-\theta)}),$ we
obtain that
$$
f\ast
\ol{h}(e^{i\theta})=\f1{2\pi}\int_0^{2\pi}f(e^{i(s-\theta)})\ol{h(e^{is})}ds=0,\forall
f\in H_1\cup H^0\cup H^-(\TT).$$ Now by Theorem \ref{con0}, we
obtain $h\equiv0$, which contradicts with $H_2$ being proper. So
$Z(H_1)\cap \ZZ^+$ can not be empty.

 According to the just proved \bf Claim \rm  there exists at least one positive integer $m$
such that $m\in Z(H_1)$. Since $g(t)=t^m\in
\tilde{H}^+(\TT),$ there exist $f'(t)\in H_1$ and $f''(t)\in
H_2$
 such that $t^m=f'(t)+f''(t)$. Recall that $m\in Z(H_1)$ means $f'_m=0$. Then we have
 \begin{equation}\label{}
   \nonumber 0=<f',f''>=<f',t^m-f'>=<f',t^m>-<f'.f'>=f'_m-||f'||^2=-||f'||^2,
 \end{equation}
 which gives $f'(t)=0$. So $t^m=f''(t)\in H_2$.  Now by (\ref{nointeger}) we have $t=\pi_{(\f1m,0)}t^m\in H_2$.
 As consequence, $t^n=\pi_{(n,0)}t$ are all in $H_2$ for  $n\in\ZZ^+$.
  Hence $\tilde{H}^+(\TT)=\ol{\mathrm{span}(t,t^2,\cdots,)}\subseteq H_2 \subseteq \tilde{H}^+(\TT)$, which contradicts with the assumption that $H_2$ is proper.

\epp

By the Plemelj formula (\cite{miklin}) it is easy to prove the following theorem.
 \bT\label{cirl} Let $\tilde{ \mathbf{H}}$  be the Hilbert transformation on $L^2(\TT)$. Then
 \begin{equation}\label{cauchy0}
  \tilde{ \mathbf{H}}|_{\tilde{ H}^+}=-i\mathbf{I}|_{\tilde{ H}^+},\ \mathbf{C}|_{H^0}=0,\ \mathbf{C}|_{H^-}=i\mathbf{I}|_{H^-}.
 \end{equation}
 \eT

 \noindent We end this paper with the inverse of Theorem \ref{cirl}.

 \bT\label{endd}
Let $T$ be a bounded operator from $L^2(\TT)$ to itself. Assume that $\tilde{ \mathbf{T}}$ commutes with $\pi_{(\a,\b)}, (\a,\b)\in G$.
 Then there exist three complex numbers $\lambda,\eta,\om$ such that
 $$
  \tilde{ \mathbf{T}}|_{\tilde{ H}^+}=\lambda \mathbf{I}|_{\tilde{ H}^+},\ \tilde{ \mathbf{T}}|_{H^0}=\eta \mathbf{I}_{H_0},\ \tilde{ \mathbf{T}}|_{H^-}=\om \mathbf{I}|_{H^-}.
 $$
 \eT
\bpp It is obvious that  the family of $\pi_{(\a,\b)},\ (\a,\b)\in G)$ acts irreducibly on $H^0$. By Theorem \ref{dim} and \ref{schur},
$T$ must be scalar operators on, respectively, $H^+,H^0$, and $H^-$.

\epp
\brs
I. Theorem \ref{endd}, as far as what we are aware, is a new result. It is also not easy to prove by the methods from \cite{king,stein} because of  lack of the general dilation on the unit circle.

 II. It is also possible to characterize the circular Hilbert transformation by the symmetry of  the generalized M\"{o}bius group containing both the rotations and M\"{o}bius transforms. Define $\tau_\theta z=e^{i\theta}z$, $\theta\in\RR$ and $\varphi_a(z)=\f{z-a}{1-az}$, $a\in(0,1)$. Then the M\"{o}bius group, $\mathcal{M},$ is generated by $\tau_\theta$ and $\varphi_a$. There exists a natural representation of $\mathcal{M}$ over $L^2(\TT)$ as following: for $\varphi\in \mathcal{M}$ with the expression $\varphi(z)=e^{i\theta}\f{z-a}{1-az}$, we define
 \begin{equation}\label{}
   \nonumber ((\pi_\mathcal{M}\varphi)(f))(t)=\f{\sqrt{1-a^2}}{1-at}f(\varphi^{-1}(t)),\ f\in L^2{(\TT)}.
 \end{equation}
It is easy to check that $\pi_\mathcal{M}$ is a unitary representation and commutates with  $\mathbf{C}$. We can also prove that its restrictions to $H^\pm$ are respectively irreducible. Then we obtain Theorem \ref{cirl}.

But the Fourier series expansions of the M\"{o}bius transforms are  complicated. So we can not obtain the precise structures of $\mathbf{C}$ in the phase space as we discussed in the real line case.

\ers

\noindent\bf Acknowledgements \rm The first author is supported by Macao Government Science and Technology Development Fund, MSAR. Ref. 045/2015/A2 and NSFC11701597. The second author is supported by NSFC11601390 and NSFC11471250. The third author is supported by Macao Government Science and Technology Development Fund FDCT099 and FDCT079.

\end{document}